\documentclass[11pt,reqno]{amsart}
\usepackage{amssymb,accents,covington}

\def\noi{\noindent}

\newtheorem{Thm}{Theorem}[section]
\newtheorem{Def}[Thm]{Definition}
\newtheorem{Lm}[Thm]{Lemma}
\newtheorem{Prop}[Thm]{Proposition}

\newtheorem{Rem}[Thm]{Remark}

\def\cal{\mathcal}
\def\Bbb{\mathbb}
\def\mf{\mathfrak}

\def\<{\left<}
\def\>{\right>}

\def\a{\alpha}
\def\b{\beta}
\def\d{\delta}

\def\th{\theta}

\def\l{\lambda}
\def\L{\Lambda}

\def\Re{\Bbb R}
\def\F{\Bbb F}
\def\C{\Bbb C}
\def\Z{\Bbb Z}

\def\H{\cal H}
\def\A{\cal A}

\def\G{\mf h}

\def\W{\ring W}
\def\w{\ring w}
\def\RR{\ring R}
\def\Q{\ring Q}

\def\Gc{\ring {\mf h}}

\begin{document}
\title[double affine Hecke algebras]
{Involutions of double affine Hecke algebras}
\author{Bogdan Ion}
\address{
Department of Mathematics,
Princeton University,
Princeton NJ-08544
}
\email{bogdan@math.princeton.edu}
\begin{abstract}
The main aim of the paper is to formulate and prove a result about the
structure of double affine Hecke algebras  
which allows  its two commutative subalgebras to play a symmetric role.
This result is essential for the 
theory of intertwiners of double affine Hecke algebras.
\end{abstract}
\maketitle


\thispagestyle{empty}
\section*{Introduction}
\bigskip
Double affine Hecke algebras were introduced by Cherednik
\cite{cherednik}, \cite{cherednik2} in connection with affine quantum
Knizhnik-Zamolodchikov equations and eigenvalue problems of Macdonald type. 
As a general principle one can associate to any type of root system
(finite, affine, elliptic) a Weyl group, an Artin group 
and a Hecke algebra. Along these 
lines, double affine Hecke algebras are associated with a certain class
of elliptic root systems. The Section \ref{preliminaries} contains the 
definitions and relevant results about double affine Weyl groups, their
Artin groups and their Hecke algebras as well as a topological interpretation
of affine Artin  groups due to H. van der Lek.

One of the key features of double affine Hecke algebras 
is that they contain two commutative 
subalgebras. Our main result roughly states that there exists a pairing 
between the root systems in question such that 
 their corresponding
Weyl groups, Artin groups and Hecke algebras are canonically isomorphic
in such a way that 
for the Hecke algebras
the roles of the above mentioned commutative subalgebras
is interchanged. The Section \ref{involution} contains the precise statements 
of these results. 
As a consequence of the existence of such canonical 
isomorphisms there exists the possibility of defining some special elements
of the double affine Hecke algebras, called { intertwiners}.
The importance of the intertwiners was first 
understood by Knop and Sahi \cite{knop}, \cite{ks}, 
\cite{s1} for ${\rm GL}_n$ and then by Cherednik \cite{cherednik3} in the 
general (reduced) case.
Their properties, as they follow from the isomorphism allows them to play
a central role in the theory of (non)symmetric Macdonald polynomials (see
\cite{cherednik3} for a full account on  intertwiners).

The proof rely on the topological interpretation of extended Artin groups,
due to H. van der Lek \cite{lek}, \cite{lek-thesis}. This relates the
double affine Artin groups and the fundamental groups of certain
complex hyperplane complements. 
Section \ref{topological} reveals this connection and establishes
the canonical isomorphism at the level of Artin groups. Section \ref{Hecke}
finishes the proof of  the existence of the canonical 
isomorphism at the level of Hecke algebras.

The statements of our main result as well as the connection between the double
affine Artin groups and certain fundamental groups
were announced in 
Theorem 2.2 and Theorem 2.4 of \cite{cherednik}.
 
{\sl Acknowledgement.} I would like  to thank Professor Siddhartha Sahi 
for many enlightening discussions on the subject and much useful advice.

\bigskip
\section{Preliminaries}\label{preliminaries}
\bigskip

\subsection{Affine and double affine Weyl groups}\label{notations}

For the most part we shall adhere to the notation in \cite{kac}.
Let $A=(a_{jk})_{0\leq j,k\leq n}$ be an irreducible \emph{affine}
Cartan matrix,
$S(A)$ the Dynkin diagram and  $(a_0,\dots, a_n)$  the numerical labels of
$S(A)$ in Table Aff from \cite{kac}, p.48-49. We denote by
$(a_0^\vee,\dots, a_n^\vee)$
the labels of the Dynkin diagram $S(A^t)$ of the
dual algebra which is obtained from $S(A)$ by reversing the direction
of all arrows and keeping the same enumeration of the vertices.
Let $({\G}, R, R^{\vee})$ be a realization of $A$
and let
$(\Gc, \RR,
\RR^{\vee})$
be the associated finite root system (which is a realization of the Cartan
matrix $\ring A = (a_{jk})_{1\leq j,k\leq n}$).
If we denote by $\{\a_j\}_{0\leq j\leq n}$ a basis of $R$
such that $\{\a_j\}_{1\leq j\leq n}$ is a basis of $\RR$
we have the following description
$$
{\G}^*={\Gc^*} + {\mathbb R}\delta +
{\mathbb R}{\Lambda}_0\ ,
$$
where $\d=\sum_{j=0}^n a_j\a_j$. The vector space ${\G}^*$ has a canonical
scalar product defined as follows
$$
(\a_j,\a_k):=d_j^{-1}a_{jk}\ ,\ \ \ \
(\L_0,\a_j):=\d_{j,0}a_0^{-1}\ \ \ \text{and}\ \ \ (\L_0,\L_0):=0,
$$
with $d_j:= a_ja_j^{{\vee}-1}$ and $\d_{j,0}$  Kronecker's delta.
As usual, $\{\a_j^\vee:=d_j\a_j\}_{0\leq j\leq n}$ are the coroots.
Denote by
$\Q =\oplus_{j=1}^n\Z\a_j$ and $\Q^\vee=\oplus_{j=1}^n\Z\a_j^\vee$
the root lattice, respectively the coroot lattice of $\RR$ and by
$Q=\oplus_{j=0}^n\Z\a_j=\Q\oplus \Z\d$ the root lattice of $R$.

\medskip
Given $\a\in R$, $x\in \G^*$ let
$$
s_\a(x):=x-\frac{2(x,\a)}{(\a,\a)}\a\ .
$$
The {\sl affine Weyl group} $W$ is the subgroup of ${\rm GL}(\G^*)$
generated by all $s_\a$ (the simple
reflexions $s_j=s_{\a_j}$ are enough).
The {\sl finite Weyl group} $\W$ is
the subgroup generated by $s_1,\dots,s_n$.
Both the finite and the affine Weyl group are
Coxeter groups and they can be abstractly
defined as generated by $s_1,\dots,s_n$, respectively
$s_0,\dots,s_n$, and some relations. These relations are called
Coxeter relations
and they are of two types:
the reflection relations $s_j^2=1$ and the braid relations
(see \cite{humph} for details).

\medskip
The {\sl double affine Weyl group} $\tilde W$
is defined to be the semidirect product $W\ltimes Q$
of the affine Weyl group 
and the lattice $Q$
(regarded as an abelian group with elements
$\tau_\b$, where $\b$ is a root), the affine
Weyl group (with the opposed multiplication)
acting on the root lattice as follows
$$
w\tau_\b w^{-1}=\tau_{w(\b)}.
$$
This group is the hyperbolic extension of an elliptic Weyl group,
that is the Weyl group associated with an elliptic root system
(see \cite{saito2} for definitions). It also has a presentation
with generators and relations (elliptic
Coxeter relations). We refer the reader to \cite{saito2}
for the details.

\medskip
The affine Weyl group $W$ can also be presented
as a semidirect product in the following way.
For $1\leq j\leq n$, let $e_j=\text{max}\{a_0^{-1}, d_j\}$.
Recall that $a_0=1$ in all cases except for
$A=A_{2n}^{(2)}$, where $a_0=2$. Denote by $M$ the lattice generated by
$\{A_j=e_j\a_j\}_{1\leq j\leq n}$. Then $W$ is
the semidirect product
of $\W$ and the lattice $M$ (regarded as an abelian group with elements
$\l_\mu$, where $\mu$ is in $M$),
the finite Weyl group acting on the lattice
$M$ as follows
$$
w\l_\mu w^{-1}=\l_{w(\mu)}.
$$
Let us remark that the numbers $e_j$ depend only on the length of the
corresponding simple root. Therefore, we will write $e_s$, $e_l$ for
$e_j$ if the corresponding simple root is short, respectively long.
The number $p:=e_s$ 
will play an important role. It is easy to see
that $p\in\{1,2,3\}$.
If $\b\in \RR$, $e_\b$ will denote $e_s$ or $e_l$ depending on the length
of $\b$. With this notation define, 
$$
A_\b:=e_\b\b \ \ \text{and}\ \ 
A_\b^\vee := e_\b^{-1}\b^\vee
.
$$
\medskip
For $r$ a real number, $\G^*_r=\{ x\in\G\ ;\ (x,\d)=r\}$ is the
level $r$ of $\G^*$. We have
$$
\G^*_r=\G^*_0+r\L_0=\Gc^*+{\Bbb R}\d+r\L_0\ .
$$
The action of $W$
preserves each of the $\G^*_r$ and  we can identify each of the $\G^*_r$
canonically with $\G^*_0$ and obtain an (affine) action of $W$ on $\G^*_0$.
For example, the level zero action of $s_0$ and $\l_\mu$ on $\G^*_0$ is
\begin{eqnarray*}
s_0(x)     & = & s_\th(x)+(x,\th)\d\ ,\\
\l_\mu(x)   & = & x - (x,\mu)\d\ ,
\end{eqnarray*}
and the level one affine action of the same elements on $\Gc^*$ is 
\begin{eqnarray*}
s_0(x)   & = & s_\th(x)+a_0^{-1}\th\ ,\\
\l_\mu(x) & = & x +\mu \ ,
\end{eqnarray*}
where we denoted by
$\th=\delta -a_0\a_0$.

\bigskip
\subsection{Artin groups and Hecke algebras}
\bigskip
To any Coxeter group we can associate its Artin group, as the group defined
with the same generators which satisfy only the  {\sl braid relations} 
(that is, forgetting the
reflexion relations). The finite and affine Weyl groups are Coxeter groups; 
we will
make precise the definition of the Artin groups in these cases. First let us
make the following convention

\bigskip
\begin{quotation}\label{restriction}
{\sl
For
the rest of the paper we assume our irreducible affine Cartan matrix to
be of any type except  $A_{2n}^{(2)}$.
}
\end{quotation}
\bigskip

\noi Note that the
corresponding definitions and results of this paper (at the level of Hecke 
algebras) 
for the case omitted 
can be found in 
\cite{sahi}.
\begin{Def} With the notation above define
\begin{enumerate}
\item[i)]
the finite Artin group $\A_{\W}$ as the group
generated by elements
$$T_1,\dots,T_n$$ 
satisfying the same braid relations as the reflexions $s_1,\dots,s_n$;
\item[ii)]
the affine Artin group $\A_W$ as the group
generated by the elements
$$T_0,\dots,T_n$$ 
satisfying the same braid relations as the reflexions $s_0,\dots,s_n$.
\end{enumerate}
\end{Def}
From the definition it is 
clear that the finite Artin group can be realized as a subgroup inside
the affine Artin group. To
define the Hecke algebras,
we introduce a field $\F$ (of parameters) as follows:
fix indeterminates $q$ and $t_0,\dots,t_n$ such that
$
t_j=t_k \text{ iff } d_j=d_k;
$
let $m$ be the lowest common denominator of the rational numbers
$\{(\a_j,\l_k)\ |\ 1\leq j,k\leq n \}$, and let $\F$ denote the field of
rational functions in $q^{1/m}$ and $t_j^{1/2}$. Because in our case there
are at most two different root lengths we will also use the notation
$t_l$, $t_s$ for $t_j$ if the corresponding simple root is long,
respectively short.

For further use we also introduce the following lattices
$\cal{M}_Y:=\{Y_\mu;\mu\in M\}$ and $\cal{Q}_X:=\{X_\b;\b\in \Q\}$.
We will use the same notation for their
group $\F$-algebras.
\begin{Def}
The finite Hecke algebra $\H_{\W}$ is the quotient of the group
$\F$-algebra of the finite Artin group by the relations
\begin{equation}
T_j-T_j^{-1}=t_j^{1/2} -t_j^{-1/2},
\end{equation}
for $1\leq j\leq n$.
\end{Def}
Recall that the affine Weyl group also has a second presentation, as a
semidirect product. There is a corresponding description of the
affine Artin group (and consequently of the
affine Hecke algebra) due to van der Lek \cite{lek} and
Lusztig \cite{lusztig}.
\begin{Prop}\label{prop2.2}
The affine Artin group $\A_W$ is generated by the finite Artin group
and the lattice $\cal{M}_Y$ such
that the following relations are satisfied for all $1\leq j\leq n$
\begin{enumerate}
\item[(i)] $T_jY_\mu=Y_\mu T_j \ \ \text{if \ } (\mu,A_j^\vee)=0$,
\item[(ii)]
$T_jY_\mu T_j=Y_{s_j(\mu)}   \ \ \text{if \ } (\mu,A_j^\vee)=1$.
\end{enumerate}
\end{Prop}
\begin{Rem}\label{antidom}
In this description $Y_\mu=T_{\l_\mu}$ for $\mu$ any
anti-dominant element of the lattice $M$. For example $Y_{-\th}=T_{s_\th}T_0$.
\end{Rem}
\begin{Def}
The affine Hecke algebra $\H_W$ is the quotient of the group
$\F$-algebra of the affine Artin group by the relations
\begin{equation}\label{q1}
T_j-T_j^{-1}=t_j^{1/2} -t_j^{-1/2},\ \ \ \text{for all}\ \  0\leq j\leq n,
\end{equation}
\begin{equation}\label{q2}
Y_{-A_j}T_j^{-1}-T_jY_{A_j}=t_j^{1/2} -t_j^{-1/2},\ \ \ 
\text{for all}\ \  1\leq j\leq n.
\end{equation}
\end{Def}
The elements $T_1, \dots, T_n$ generate the  {\sl finite Hecke algebra}
$\H_{\W}$. There are natural bases of $\H_W$ and $\H_{\W}$:
$\{T_w\}_w$ indexed by
$w$ in $W$ and in
$\W$ respectively, where $T_w=T_{j_l}\cdots T_{j_1}$ if
$w=s_{j_l}\cdots s_{j_1}$ is a reduced  expression of $w$
in terms of simple reflexions.

\medskip

We must remark that 
the objects we define here are not the objects traditionally used
in the literature
:  {\sl extended} affine Artin groups and Hecke algebras.
Nevertheless, their structure is completely similar.
Regarding the definition of Hecke algebras we note that
the relations (\ref{q2}) and (\ref{q1})
are conjugate inside the extended
affine Artin group, but not inside the affine Artin group. Therefore,
only the relations (\ref{q1}) have to be required in the definition of
extended affine Hecke algebras.

\medskip

As before, the affine Hecke algebra has a presentation in terms of the
finite Hecke algebra and the group algebra $\cal{M}_Y$ (see \cite{lusztig}
for the
proof).
\begin{Prop}\label{gen-lusztig}
The affine Hecke algebra $\H_W$ is generated by the finite Hecke algebra
and the group algebra $\cal{M}_Y$ such that the following relation is
satisfied
for any
$\mu$ in the lattice $M$ and
any $1\leq j\leq n$\ :
$$
Y_\mu T_j-T_jY_{s_j(\mu)} = (t_j^{1/2}-t_j^{-1/2})
\frac{Y_\mu-Y_{s_j(\mu)}}{1-Y_{A_j}}\ .
$$
\end{Prop}
\begin{Rem} An immediate consequence of the Proposition \ref{gen-lusztig}
is that the relation
$$
T_0-T_0^{-1}=t_0^{1/2}-t_0^{-1/2}
$$
it is contained in the ideal generated by all the other relations defining 
the affine Hecke algebra.
\end{Rem}
The definition of the double affine Artin group and 
Hecke algebra is due
to Cherednik (see, for example, \cite{cherednik}).
\begin{Def}
The double affine Artin group $\A_{\tilde W}$ is
generated by
the affine Artin group $\A_W$, the lattice $\cal{Q}_X$ and the element
$X_\d$ such the following relations are satisfied for all
$0\leq i\leq n$
\begin{enumerate}
\item[(i)] $X_\d$ is central\ ,
\item[(ii)] $T_jX_\b=X_\b T_j \ \ \text{if \ } (\b,\a_j^\vee)=0$,
\item[(iii)] $T_jX_\b T_j=X_{s_j(\b)}   \ \ \text{if \ } (\b,\a_j^\vee)=-1.$
\end{enumerate}
\end{Def}

The double affine Weyl group it is not a Coxeter group, but a generalized
Coxeter group (in the sense of Saito and
Takebayashi, see \cite{saito2}) and we can define
the associated Artin group in the same way as for a Coxeter group (that is,
by keeping the generalized braid relations and forgetting the reflexion
relations). 
By our knowledge
at present time the literature does not contain
any result that establishes the equivalence of the two definitions.
\begin{Def}\label{def3}
The double affine Hecke algebra $\H_{\tilde W}$ is the quotient of the group
$\F$-algebra of the double affine Artin group by the relations
\begin{equation}\label{def1}
T_j-T_j^{-1}=t_j^{1/2} -t_j^{-1/2},\ \ \ \text{for all } \ \ 1\leq j\leq n,
\end{equation}
\begin{equation}\label{def1.2}
Y_{-A_j}T_j^{-1}-
T_jY_{A_j}=t_j^{1/2} -t_j^{-1/2},\ \ \ \text{for all } \ \ 1\leq j\leq n,
\end{equation}
\begin{equation}\label{def1.3}
T_j^{-1}X_{\a_j}-X_{-\a_j}T_j=t_j^{1/2} -t_j^{-1/2},
\ \ \ \text{for all } \ \ 1\leq j\leq n,
\end{equation}
\begin{equation}\label{eq7}
T_0^{-1}X_{\a_0}-X_{-\a_0}T_0=t_0^{1/2} -t_0^{-1/2},
\end{equation}
 and by 
\begin{equation}
X_\d=q^{-1}.
\end{equation}
\end{Def}
\noi As before, $\H_{\tilde W}$ can be described in terms of $\H_W$ and
the group algebra $\cal{Q}_X$.
\begin{Prop}[\cite{cherednik}]
The double affine Hecke algebra
$\H_{\tilde W}$ is the $\F$-algebra generated by
the affine Hecke algebra $\H_W$ and the group algebra $\cal{Q}_X$
such that, with the notation $X_\d=q^{-1}$, 
the following relation is satisfied for any
$\b$ in the root lattice and
any $0\leq j\leq n$:
\begin{equation}\label{def2}
T_jX_\b-X_{s_j(\b)}T_j = (t_j^{1/2}-t_j^{-1/2})
\frac{X_\b-X_{s_j(\b)}}{1-X_{-\a_j}}\ .
\end{equation}
\end{Prop}
The proof is completely similar with the proof of the Proposition 
\ref{gen-lusztig}.
\subsection{Some results of van der Lek}

Inspired by the techniques introduced in \cite{deligne},
H. van der Lek developed 
in \cite{lek-thesis}
a machinery which allows one 
to compute fundamental groups of complex hyperplane complements. As an 
application he realized the affine Artin groups as fundamental groups
of such spaces. Below, we will briefly review his result.
 
With the notation from Section \ref{notations} let
$$V=\Gc^*.$$
The
finite Weyl group acts properly discontinuously on $V$ and in
consequence the diagonal action of $\W$ on the domain
$
V+iV
$
has the same property (here $i=\sqrt{-1}$).

For $\b\in \RR$ any root 
and any integer $k$ denote by
$$
\bar H_{A_{\b},k}= \{  v\in  V\ |
\ (v,A_{\b})=k(A_{\b},A_{\b})\ \}.
$$
Consider the domain
$$
\bar \Omega=V+iV 
$$
and the following action of $W$ on it:
\begin{enumerate}
\item[$\bullet$] $s_j(v_1+iv_2)=
s_j(v_1)+is_j(v_2)$ for $j\neq 0$ ;
\item[$\bullet$]
$\l_\mu(v_1+iv_2)=
v_1 +\mu+iv_2$\ 
for $\mu\in M$.
\end{enumerate}
 
Consider the following space 
$$
\bar{\cal{Y}}:=\left(\left(V+iV\right)\ \ -
\bigcup_{\stackrel{\b\in \RR,}{\scriptscriptstyle{k\in\Z}}}
\left(\bar H_{A_{\b},k}+i\bar H_{A_{\b},0} \right)\right).
$$
and the orbits space associated to the above action 
$$
\bar{\cal{X}}=\bar{\cal{Y}}/{W}. 
$$
Let $\bar p:\bar{\cal{Y}}\to \bar{\cal{X}}$ denote the canonical
projection. Fix  $\bar c\in V$ such that 
$\bar c$ is in the fundamental chamber for the action of $\W$.
We can choose
$$
\bar \star=\bar p(i\bar c)
$$ 
as a base point for
$\bar{\cal{X}}$.
To state Theorem 2.5, Section 3 of 
\cite{lek-thesis} we need the following
notation 
\begin{enumerate}
\item[$\bullet$] $\bar \Upsilon_j:[0,1]\to \bar{\cal{Y}}$;
$\bar \Upsilon_j(t)=\bar c +\frac{e^{-\pi it}-1}{2}(\bar c,\a_j^\vee)\a_j$, 
$1\leq j\leq n$;
\item[$\bullet$] $\bar y_{A_j}:[0,1]\to \bar{\cal{Y}}$ ;
$\bar y_{A_j}(t)=\bar c+tA_j$, $1\leq j\leq n$.
\end{enumerate}
Note that, $\bar p\circ \bar \Upsilon_j$ and
$\bar p\circ \bar y_{A_j}$  
are closed paths
in $\bar{\cal{X}}$. 
Now, Theorem 2.5, Section 3 of
\cite{lek-thesis} states as follows.
\begin{Thm}
With the notation above, the fundamental group 
$\pi_1(\bar{\cal{X}},\bar \star)$
and the affine Artin group $\A_{W}$ are isomorphic.
Under the isomorphism the homotopy classes of the paths 
$\bar p\circ \bar \Upsilon_j$, 
and $\bar p\circ\bar y_{A_j}$
correspond to $T_j$ and $Y_{A_j}$, respectively.
\end{Thm}
In fact van der Lek's Theorem 2.5, Section 3 is more general. It can be stated 
for an arbitrary finite Coxeter group. Although the statements don't include
the eventual existence of imaginary roots in the root system associated with
this Coxeter group, once one appropriately modifies the regular orbits space
the proofs work exactly the same way to produce a corresponding result. 
In Section \ref{topological} we will
define the space which has the fundamental group isomorphic with 
the double affine Artin group.

\bigskip
\section{The involution}\label{involution}
\bigskip
In this section we will define an involution 
$\iota$ on the set with elements of the form $(A, \{\a_j\})$
where $A$ is a 
irreducible affine Cartan matrix satisfying our
convention in Section \ref{restriction} and $\{\a_j\}$ is a basis of the 
corresponding root system.

Because the
irreducible root systems are classified by their Dynkin diagrams it is
enough to make precise a basis of $R^{\iota}$, which we denote by
$\{\a_j^{\iota}\}_{0\leq j\leq n}$
and the scalar products between the simple
roots.
The fact that we want our involution to
preserve the affine type allows us to specify only
$\{\a_j^{\iota}\}_{1\leq j\leq n}$. For $1\leq j\leq n$, define
$$
\a_j^{\iota}=\frac{A_j}{\sqrt{p}}\ ,
$$
the scalar product being the canonical one in $\G^*$. It will cause
no confusion we will realize the root system $R^{\iota}$ on $\G^*$. Also,
$$
\d^{\iota}=\d \ \ \text{and}\ \ \L_0^{\iota}=\L_0.
$$

The finite Weyl groups associated to $\RR$ and $\RR^\iota$ coincide and we
will not distinguish between them.
The affine and double affine Weyl groups
associated with $R^{\iota}$ we denote by $W^{\iota}$,
respectively $\tilde W^{\iota}$.

At the level of affine root systems the involution fixes all irreducible
affine Cartan matrices, except $B_n^{(1)}$ and $C_n^{(1)}$ which are 
interchanged. At the level of roots the involution acts as 
identity for the all root systems except $B_n^{(1)}$, $C_n^{(1)}$, 
$F_4^{(1)}$ and $G_2^{(1)}$. For root systems of these types
these the involution interchanges short and 
long roots. 
At the level of double affine Weyl groups the following result holds.
\begin{Prop}
With the notation above $M^{\iota}=\sqrt{p}\Q$ and
$\Q^{\iota}=\frac{1}{\sqrt{p}}M$. Furthermore, the map
$$
\phi_{\tilde W}:\tilde W\to \tilde W^{\iota}
$$
defined as follows
\begin{eqnarray*}
\phi_{\tilde W}(\w)            & = &  \w\ \ \ \text{ for } \w\in \W\ ,\\
\phi_{\tilde W}(\l_\mu)       & = &  \tau^{\iota}_{\frac{1}{\sqrt{p}}\mu}
\ \ \ \text{ for }\mu\in M\ ,\\
\phi_{\tilde W}(\tau_\beta)   & = &  \l^{\iota}_{\sqrt{p}\b}
\ \ \ \text{ for } \b\in \Q\ ,\\
\phi_{\tilde W}(\tau_\delta)  & = &  \tau^{\iota}_{-\d}\ .
\end{eqnarray*}
is an isomorphism of groups. The inverse of $\phi_{\tilde W}$ is
$\phi_{\tilde W^{\iota}}$.
\end{Prop}
\begin{proof}
All computations are straightforward.
\end{proof}

The main results of the paper state that the isomorphism between the
double affine Weyl groups from the previous Proposition induces
isomorphisms
at the level of double affine Artin groups and double
affine Hecke algebras. In order to keep our notation as simple as possible
we denoted all these isomorphisms by the same symbol.
\begin{Thm}\label{teorema1}
The map
\begin{eqnarray*}
\phi_{\tilde W}(T_j)            & = &  {T^{\iota}_j}^{-1}\ \ \ \text{ for }
1\leq j\leq n\ ,\\
\phi_{\tilde W}(Y_\mu)       & = &  X^{\iota}_{\frac{1}{\sqrt{p}}\mu}
\ \ \ \text{ for }\mu\in M\ ,\\
\phi_{\tilde W}(X_\beta)   & = &  Y^{\iota}_{\sqrt{p}\b}
\ \ \ \text{ for } \b\in \Q\ ,\\
\phi_{\tilde W}(X_\delta)  & = &  X^{\iota}_{-\d}\ ,
\end{eqnarray*}
can be uniquely extended to an isomorphism
$$
\phi_{\tilde W}:\A_{\tilde W}\to \A_{\tilde W^{\iota}}\ ,
$$
between the double affine Artin
groups. The inverse of $\phi_{\tilde W}$ is
$\phi_{\tilde W^{\iota}}$.
\end{Thm}
\begin{Thm}\label{teorema2}
The group isomorphism
from the previous Theorem  extended $\C$-linearly and by
\begin{eqnarray*}
\phi_{\tilde W}(t_j)   & = &  t_j^{-1}\ \ \text{for}\ j\neq 0,\\
\phi_{\tilde W}(q)  & = &  q^{-1}\ .
\end{eqnarray*}
induces an isomorphism between the corresponding Hecke algebras.
The inverse of $\phi_{\tilde W}$ is
$\phi_{\tilde W^{\iota}}$.
\end{Thm}

\noi The proofs of Theorem \ref{teorema1} and Theorem \ref{teorema2} will
be concluded in
Section \ref{topological} and Section \ref{Hecke}, respectively.
\bigskip
\section{The topological interpretation}\label{topological}
\bigskip
\subsection{The orbits space}\label{orbitspace}
The proof of the Theorem \ref{teorema1}
rely on the interpretation of the double affine Artin group
$\A_{\tilde W}$ as a fundamental group of a certain topological space. 
In what follows we will present the construction of this space.  

With the notation from Section \ref{notations} let
$$
V=\Gc^* \ \ \text{and}\ \  \tilde V=\Gc^*+\Re\L_0.
$$
The Tits cone is defined to be
$$
I:=\bigcup_{w\in W}w(\bar C),
$$
where
$$
C:=\{ \tilde v\in \tilde V\ | (\tilde v,\a_j)> 0\ ,\ 0\leq j\leq n\}
$$
is the fundamental chamber for the action of the affine Weyl group $W$. 
The interior of the Tits cone is
$$
\ring I=\{v+r\L_0\ | v\in V,\ r>0\ \}.
$$
Recall that $x\in I$ is an interior point iff ${\rm stab}_W(x)$
is a finite group. The
affine Weyl group acts properly discontinuously on $\ring I$ and in
consequence the diagonal action of $W$ on the domain
$
\tilde V+i\ring I
$
has the same property.

For $\b\in R^{re}$ any real root (that is, non-proportional with the imaginary
root $\d$)  and any integer $k$ denote by
$$
\tilde H_{\b,k}= \{ \tilde v\in \tilde V\ |\ (\tilde v,\b^\vee)=k\ \}.
$$
Consider the domain
$$
\tilde\Omega=(\tilde V+i\ring I) \times\C
$$
and the following action of $\tilde W$ on it:
\begin{enumerate}
\item[$\bullet$] $s_j(\tilde v_1+i\tilde v_2;z)=
(s_j(\tilde v_1)+is_j(\tilde v_2);z)$ for $j\neq 0$ ;
\item[$\bullet$]$s_0(\tilde v_1+i\tilde v_2;z)=
(\tilde v_1+(\tilde v_1,\a_0)\th+i\tilde v_2+i(\tilde v_2,\a_0)\th;
ze^{\frac{\pi i}{2} \{ (\tilde v_1,\th) -
\frac{(\tilde v_1,\d)}{(\tilde v_2,\d)}(\tilde v_2,\th) \} } )$;
\item[$\bullet$]
$\tau_\b(\tilde v_1+i\tilde v_2;z)=
(\tilde v_1 +\b+i\tilde v_2;
ze^{\frac{\pi i}{2} \frac{(\tilde v_2,\b)}{(\tilde v_2,\d)}})$
for $\b\in \Q$;
\item[$\bullet$]
$\tau_\d(\tilde v_1+i\tilde v_2;z)=
(\tilde v_1+i\tilde v_2,-z)$.
\end{enumerate}
 
Consider the following space 
$$
\tilde{\cal{Y}}:=\left(\left(\tilde V+i\ring I\right)\ \ -
\bigcup_{\stackrel{\b\in R^{re},}{\scriptscriptstyle{k\in\Z}}}
\left(\tilde H_{\b,k}+i\tilde H_{\b,0} \right)\right)\times \C^*.
$$
and the orbits space associated to the above action is
$$
\tilde{\cal{X}}=\tilde{\cal{Y}}/{\tilde W}. 
$$
Let $\tilde p:\tilde{\cal{Y}}\to \tilde{\cal{X}}$ denote the canonical
projection. Fix  $c\in V$ such that $c+\L_0$ is in the fundamental chamber 
and the 
numbers $\{(c,\a_j)\}_{1\leq j\leq n}$
are positive and sufficiently small. Fix also $z_0$ a positive real number
and let $u=c+ic$.
Then the point 
$(u+i\L_0;z_0)$ is in $\tilde{\cal{Y}}$ and we can choose
$$
\tilde \star=\tilde p(u+i\L_0;z_0)
$$ 
as a base point for
$\tilde{\cal{X}}$. To state Theorem 2.5, Section 3 of 
\cite{lek-thesis} adapted to our present situation we need the following
notation 
\begin{enumerate}
\item[$\bullet$] $\tilde \Upsilon_j:[0,1]\to \tilde{\cal{Y}}$;
$\tilde \Upsilon_j(t)=(u+i\L_0 +\frac{e^{\pi it}-1}{2}(u,\a_j^\vee)\a_j;z_0)$, 
$1\leq j\leq n$;
\item[$\bullet$] $\tilde \Upsilon_0:[0,1]\to \tilde{\cal{Y}}$;
$\tilde \Upsilon_0(t)=(u+i\L_0 -\frac{e^{\pi it}-1}{2}(u+i\L_0,\a_0)\th;
z_0e^{\frac{\pi i t}{2} (c,\th)})$;
\item[$\bullet$] $\tilde x_{\a_j}:[0,1]\to \tilde{\cal{Y}}$ ;
$\tilde x_{\a_j}(t)=(u+i\L_0+t\a_j;z_0e^{
\frac{\pi it}{2}(c,\b)})$, $1\leq j\leq n$;
\item[$\bullet$] $\tilde x_{\d}:[0,1]\to \tilde{\cal{Y}}$ ;
$\tilde x_{\d}(t)=(u+i\L_0;z_0e^{\pi it})$.
\end{enumerate}
Note that, $\tilde p\circ \tilde \Upsilon_j$, 
$\tilde p\circ \tilde x_{\a_j}$ and 
$\tilde p\circ \tilde x_\d$ are closed paths
in $\tilde{\cal{X}}$. Now, Theorem 2.5, Section 3 of
\cite{lek-thesis} for the affine Weyl group reads as follows.
\begin{Thm}\label{t25lek}
With the notation above, the fundamental group 
$\pi_1(\tilde{\cal{X}},\tilde \star)$
and the double affine Artin group $\A_{\tilde W}$ are isomorphic.
Under the isomorphism the homotopy classes of the paths 
$\tilde p\circ \tilde \Upsilon_j$, 
$\tilde p\circ\tilde x_{\a_j}$ and $\tilde p\circ\tilde x_\d$
correspond to $T_j$, $X_{\a_j}$ and $X_\d$, respectively.
\end{Thm}

What will be crucial in the proof of the Theorem \ref{teorema1} is another 
presentation of the double affine Artin group as a fundamental group of a 
slightly different topological space. Consider the subspace
$$
\Omega=(V+iV)\times \C
\hookrightarrow \tilde\Omega\ \ ;\ \ (v_1+iv_2;z)\mapsto (v_1+iv_2+i\L_0;z)\ .
$$
It is a simple fact that 
$\Omega$ is invariant under the action of
$\tilde W$. Let us make this action explicit:
\begin{enumerate}
\item[$\bullet$] $s_j(v_1+iv_2;z)=
(s_j(v_1)+is_j(v_2);z)$ for $j\neq 0$ ;
\item[$\bullet$]$s_0(v_1+iv_2;z)=
(s_\th(v_1)+is_\th(v_2)+i\th;
ze^{\frac{\pi i}{2}(v_1,\th)})$;
\item[$\bullet$]
$\l_\mu(v_1+iv_2;z)=
(v_1+iv_2+i\mu;ze^{-\frac{\pi i}{2}(v_1,\mu)})$ for $\mu\in M$;
\item[$\bullet$]
$\tau_\b(v_1+iv_2;z)=
(v_1 +\b+iv_2;ze^{\frac{\pi i}{2}(v_2,\b)})$
for $\b\in \Q$;
\item[$\bullet$]
$\tau_\d(v_1+i v_2;z)=
( v_1+i v_2,-z)$.
\end{enumerate}
For any nonzero $\b\in V$ and $k\in \Z$ define the hyperplanes
$$
H_{\b,k}:= \{ v\in  V\ |\ \frac{2( v,\b)}{(\b,\b)}=k\ \}.
$$
As before, we can consider the space
$$
\cal{Y}:=\left(\left(V+iV\right)\ \ -
\bigcup_{\stackrel{\b\in \RR,}{\scriptscriptstyle{h,k\in\Z}}}
\left( H_{\b,h}+i H_{A_\b,k} \right)\right)\times \C^*,
$$
the orbits space
$\cal{X}=\cal{Y}/{\tilde W}$ and the canonical projection
$p:\cal{Y}\to \cal{X}$. 
The point 
$(u;z_0)$ is in ${\cal{Y}}$ and we can choose
$$
\star=p(u;z_0)
$$ 
as a base point for
${\cal{X}}$.
  
We can define a map 
$$
ret: \tilde{\cal{Y}}\to \cal{Y},
$$
as follows
$$
ret(\tilde v_1+i\tilde v_2;z)=\left(\tilde v_1
-\frac{(\tilde v_1,\d)}{(\tilde v_2,\d)}\tilde v_2
+i\frac{1}{(\tilde v_2,\d)}\tilde v_2-i\L_0\ ;\ 
z\right).
$$
Simple computations show that this map is well defined.
\begin{Prop}\label{ret}
The map
$$
ret: \tilde{\cal{Y}}\to \cal{Y},
$$
defined above
is a deformation retract. Moreover, $ret$ is $\tilde W$ equivariant
and consequently $\tilde{\cal{X}}$ and $\cal{X}$ have the same homotopy type.
\end{Prop}
\begin{proof}
All computations are straightforward.
\end{proof}
Now, the Theorem \ref{t25lek} takes a more symmetric form. Let us define first
the ingredients:
\begin{enumerate}
\item[$\bullet$] $\Upsilon_j:[0,1]\to \cal{Y}$;
$\Upsilon_j(t)=(u+\frac{e^{\pi it}-1}{2}(u,\a_j^\vee)\a_j;z_0)$,
$1\leq j\leq n$;
\item[$\bullet$] $\Upsilon_0:[0,1]\to \cal{Y}$;
$\Upsilon_0(t)=(u-\frac{e^{\pi it}-1}{2}(u+i\L_0,\a_0)\th;
z_0e^{\frac{\pi i t}{2} (c,\th)})$;
\item[$\bullet$] $y_{A_j}:[0,1]\to \cal{Y}$ ;
$y_{A_j}(t)=(u+itA_j;z_0e^{-\frac{\pi it}{2}(c,A_j)})$, $1\leq j\leq n$;
\item[$\bullet$] $x_{\a_j}:[0,1]\to \cal{Y}$ ;
$x_{\a_j}(t)=(u+t\a_j;z_0e^{\frac{\pi it}{2}(c,\a_j)})$, $1\leq j\leq n$;
\item[$\bullet$] $x_{\d}:[0,1]\to {\cal{Y}}$ ;
$x_{\d}(t)=(u;z_0e^{\pi it})$.
\end{enumerate}
Note that, $p\circ \Upsilon_j$,  $p\circ y_{A_j}$, $p\circ x_{\a_j}$ and 
$p\circ  x_\d$ are closed paths
in ${\cal{X}}$.
\begin{Thm}\label{t25leksymm}
With the notation above, the fundamental group 
$\pi_1({\cal{X}},\star)$
and the double affine Artin group $\A_{\tilde W}$ are isomorphic.
Under the isomorphism the homotopy classes of the paths 
$p\circ \Upsilon_j$, $p\circ y_{A_j}$, 
$p\circ x_{\a_j}$ and $p\circ x_\d$
correspond to $T_j$, $Y_{A_j}$, $X_{\a_j}$ and $X_\d$, respectively.
\end{Thm}
\begin{proof}
It follows from Theorem \ref{t25lek} and Proposition \ref{ret} that the 
fundamental group $\pi_1({\cal{X}},\star)$
and the double affine Artin group $\A_{\tilde W}$ are isomorphic.
Under the isomorphism the homotopy classes of the paths 
$p\circ ret(\tilde \Upsilon_{j})=p\circ \Upsilon_{j}$ (for $0\leq j\leq n$),
$p\circ ret(\tilde x_{\a_j})=p\circ x_{\a_j}$ (for $1\leq j\leq n$) and
$p\circ ret(\tilde x_{\d})=p\circ x_{\d}$ correspond to $T_j$, $X_{\a_j}$ and
$X_\d$, respectively. Moreover, by a computation  
very
similar to the one done in the proof of Theorem 5.5,
Section 3 in \cite{lek-thesis} one can show that 
this is also an isomorphism between the subgroup of 
$\pi_1(\cal{X},\star)$ generated by 
$p\circ y_{A_j}$ (for $1\leq j\leq n$) and 
$p\circ \Upsilon_{j}$ 
(for $1\leq j\leq n$) and the affine Artin group $\A_W$ and 
that $p\circ y_{A_j}$ correspond to $Y_{A_j}$ under this isomorphism.
\end{proof}

\medskip
\subsection{Proof of Theorem \ref{teorema1}}
Let us start by an analysis of the constructions in Section \ref{orbitspace}
for the group $\tilde W^\iota$. Because 
$$
M^\iota=\sqrt{p}\Q \ \ \ \text{and}\ \ \ \Q^\iota=\frac{1}{\sqrt{p}}M
$$
we obtain that there exists a one-to-one correspondence between the finite 
root systems $\RR$ and $\RR^\iota$ which assigns to any $\b\in\RR$ a root
$\b^\iota\in \RR^\iota$ such that 
$$
\b^\iota=\frac{1}{\sqrt{p}}A_\b\ \ \ \text{and} \ \ \ A_{\b^\iota}=\sqrt{p}\b.
$$
Therefore, we have that 
$$
H_{\b^\iota,k}=\frac{1}{\sqrt{p}}H_{A_\b,k}
\ \ \ \text{and} \ \ \ 
H_{A_{\b^\iota},h}=\sqrt{p}H_{\b,h},
$$
and in consequence the topological space
$$
\cal{Y}^\iota=\left(\left(V+iV\right)\ \ -
\bigcup_{\stackrel{\b\in \RR,}{\scriptscriptstyle{h,k\in\Z}}}
\left( \frac{1}{\sqrt{p}}H_{A_\b,k}+i \sqrt{p}H_{\b,h} 
\right)\right)\times \C^*.
$$
If we denote by 
$$
\cal{Z}=\left(\left(V+iV\right)\ \ -
\bigcup_{\stackrel{\b\in \RR,}{\scriptscriptstyle{h,k\in\Z}}}
\left( H_{A_\b,k}+i H_{\b,h} 
\right)\right)\times \C^*,
$$
we see that the map 
$$
\cal{Y}^\iota\to \cal{Z},\ \ \ (v_1+iv_2;z)\mapsto 
(\sqrt{p}v_1+i\frac{1}{\sqrt{p}}v_2;z)
$$
is a homeomorphism.
We can push forward the action of $\tilde W^\iota$ obtaining the
following formulas for the action on $\cal{Z}$:
\begin{enumerate}
\item[$\bullet$] $s_j^\iota(v_1+iv_2;z)=
(s_j(v_1)+is_j(v_2);z)$ for $j\neq 0$ ;
\item[$\bullet$]
$\l_{\sqrt{p}\b}^\iota(v_1+iv_2;z)=
(v_1+iv_2+i\b;ze^{-\frac{\pi i}{2}(v_1,\b)})$ for $\b\in \Q$;
\item[$\bullet$]
$\tau_{\frac{1}{\sqrt{p}}\mu}^\iota(v_1+iv_2;z)=
(v_1 +\mu+iv_2;ze^{\frac{\pi i}{2}(v_2,\mu)})$
for $\mu\in M$;
\item[$\bullet$]
$\tau_\d^\iota(v_1+i v_2;z)=
( v_1+i v_2,-z)$.
\end{enumerate}
By Theorem \ref{t25leksymm} and the above considerations the double affine 
Artin group $\A_{\tilde W^\iota}$ is isomorphic with the fundamental group of 
$\cal{X}^\prime$
the 
orbits space associated with $\cal{Z}$ and the above action of 
$\tilde W^\iota$. To make this precise let us define the following 
paths on $\cal{Z}$:
\begin{enumerate}
\item[$\bullet$] $\Upsilon_j^\iota:[0,1]\to \cal{Z}$;
$\Upsilon_j^\iota(t)=(u+\frac{e^{\pi it}-1}{2}(u,\a_j^\vee)\a_j;z_0)$,
$1\leq j\leq n$;
\item[$\bullet$] $y_{\sqrt{p}\a_j}^\iota:[0,1]\to \cal{Z}$ ;
$y_{\sqrt{p}\a_j}^\iota(t)
=(u+it\a_j;z_0e^{-\frac{\pi it}{2}(c,\a_j)})$, $1\leq j\leq n$;
\item[$\bullet$] $x_{\frac{1}{\sqrt{p}}A_j}^\iota:[0,1]\to \cal{Z}$ ;
$x_{\frac{1}{\sqrt{p}}A_j}^\iota(t)=
(u+tA_j;z_0e^{\frac{\pi it}{2}(c,A_j)})$, $1\leq j\leq n$;
\item[$\bullet$] $x_{\d}^\iota:[0,1]\to {\cal{Z}}$ ;
$x_{\d}^\iota(t)=(u;z_0e^{\pi it})$.
\end{enumerate}
By $p^\iota:\cal{Z}\to \cal{X}^\prime$ we denote the canonical projection. The
Theorem \ref{t25leksymm} reads now as follows.
\begin{Thm}\label{t25leksymmiota}
With the notation above, the fundamental group 
$\pi_1({\cal{X}^\prime},\star)$
and the double affine Artin group $\A_{\tilde W^\iota}$ are isomorphic.
Under the isomorphism the homotopy classes of the paths 
$p^\iota\circ \Upsilon_j^\iota$, $p^\iota\circ y_{\sqrt{s}\a_j}^\iota$, 
$p^\iota\circ x_{\frac{1}{\sqrt{s}}A_j}^\iota$ and $p^\iota\circ x_\d^\iota$
correspond to $T_j^\iota$, $Y_{\sqrt{s}\a_j}^\iota$, 
$X_{\frac{1}{\sqrt{s}}A_j}^\iota$ and $X_\d^\iota$, respectively.
\end{Thm}
Looking back at the way $\cal{X}$ was defined we see that the map
$$
\cal{Y}\to\cal{Z}, \ \ \ (v_1+iv_2;z)\mapsto (v_2+iv_1;\bar z),
$$
induces an isomorphism of fundamental groups of the
associated orbits spaces (by $\bar z$ we denoted the 
complex conjugate of $z$). A straightforward analysis shows that the
homotopy classes of the paths 
$p\circ \Upsilon_j$, $p\circ y_{A_j}$, 
$p\circ x_{\a_j}$, $p\circ x_\d$ and 
$(p^\iota\circ \Upsilon_j^\iota)^{-1}$, 
$p^\iota\circ x_{\frac{1}{\sqrt{p}}A_j}^\iota$
$p^\iota\circ y_{\sqrt{p}\a_j}^\iota$,  
$(p^\iota\circ x_\d^\iota)^{-1}$ correspond respectively. 
Combining this with the Theorem \ref{t25leksymm} and Theorem 
\ref{t25leksymmiota} we finish our proof of the Theorem \ref{teorema1}.

\section{Descent to Hecke algebras}\label{Hecke}
\subsection{Some combinatorial results}\label{combinatorics}
Let us first establish some notation.
For each $w$ in $W$ let $l(w)$ be the length of a
reduced (i.e. shortest) decomposition of $w$ in terms of the $s_j$.
We have
\begin{equation}\label{length}
l(w)=|\Pi(w)|
\end{equation}
where
\begin{equation}\label{pi}
\Pi(w)=\{\a\in R_+\ |\ w(\a)\in R_-\}\ .
\end{equation}
If $w=s_{j_p}\cdots s_{j_1}$ is a reduced decomposition, then
$$
\Pi(w)=\{\a^{(k)}\ |\ 1\leq k\leq p\},
$$
with $\a^{(k)}=s_{j_1}\cdots s_{j_{k-1}}(\a_{j_k})$.

For its the basic properties of the length function on Coxeter groups see
\cite{humph}. Let us list the most important ones:
\begin{enumerate}
\item For each $0\leq j\leq n$ we have $l(s_j w)=l(w)\pm 1$ ;
\item If $l(s_j w)=l(w)-1$ then $s_jw$ can be obtained from a  certain
reduced
decomposition of $w$
by omitting a  $s_j$ factor.
\end{enumerate}
When $w\in W$ can be written as $\w\l_\mu$, with $\w\in \W$ and
$\mu \in M$ the formula (\ref{length}) takes the following form
(see \cite{lusztig})
\begin{equation}\label{length-lusztig}
l(\w\l_\mu)\ \ =\sum_{\stackrel
{\a\in \RR_+,}
{\scriptscriptstyle{\w(\a)\in \RR_-}}
}
|\frac{(\mu,\a^\vee)}{e_\a}+1|\ \  +
\sum_{\stackrel
{\a\in \RR_+,}
{\scriptscriptstyle{\w(\a)\in \RR_+}}
}
|\frac{(\mu,\a^\vee)}{e_\a}|
\end{equation}

\medskip
Next, we will consider an application of the above formula.
In this subsection we will suppose that
our irreducible affine Cartan matrix satisfies the restriction imposed
in Section \ref{restriction} and moreover
it is such that $p\neq 1$. Precisely
in this case $\th$ is the highest root of the associated finite root
system and it is not equal to $\th_s$ the highest {\sl short} root
of the associated finite root system.
\begin{Lm}\label{formula1}
Let $\A_W$ be the Artin group associated to an irreducible affine Cartan
matrix as above. Then, in $\A_W$ we have
$$
Y_{-\th^\vee_s}=T_{s_{\th_s}}T_0T_{s_{\th-\th_s}}T_0.
$$
\end{Lm}
\begin{proof}
First,
let us see that the formula in the statement makes sense. When
$p\neq 1$ all $e_\a$ equal $d_\a$, therefore $e_\a^{-1}\a^\vee=\a$ and
$M=\Q^\vee$. In consequence $-\th_s^\vee$
is an anti-dominant element of $M$.
Moreover, using basic facts about root systems, which can be found for
example in \cite[VI, \S 1, 3]{bourbaki}, we obtain that $(\th,\th_s)=1$,
which implies that $\th-\th_s=s_\th(-\th_s)\in \RR$.

Denote by $w=s_0s_{\th-\th_s}s_0$. A simple computation shows that
$s_{\th_s}w=\l_{-\th_s^\vee}$.
In the view of Remark \ref{antidom} our statement
follows from the following formulas
\begin{equation}\label{eq1}
l(w)=l(s_{\th-\th_s})+2
\end{equation}
and
\begin{equation}\label{eq2}
l(w)+l(s_{\th_s})=l(\l_{-\th_s^\vee}).
\end{equation}
The equation (\ref{eq1})
immediately follows from the second property of the
length function mentioned above keeping in mind that $w\not\in \W$,
fact which can be easily
checked. Writing down formula (\ref{length-lusztig}) for
$\l_{-\th_s^\vee}$, $s_{\th_s}$ and for 
$w=s_{\th_s}\l_{-\th_s^\vee}$ we obtain
$$
l(\l_{-\th_s^\vee})\ \ =\sum_{\stackrel
{\a\in \RR_+,}
{\scriptscriptstyle{(\a,\th_s^\vee)\neq 0}}
}
(\th_s^\vee,\a)\ ,
$$
$$
l(s_{\th_s})\ \ =\sum_{\stackrel
{\a\in \RR_+,}
{\scriptscriptstyle{(\a,\th_s^\vee)\neq 0}}
}
1 \ ,
$$
$$
l(s_{\th_s}\l_{-\th_s^\vee})\ \ =\sum_{\stackrel
{\a\in \RR_+,}
{\scriptscriptstyle{(\a,\th_s^\vee)\neq 0}}
}
\{(\th_s^\vee,\a)-1 \}\ .
$$
These formulas show that equation (\ref{eq2}) holds. The proof
statement is completed.
\end{proof}
The following fact will be crucial.
\begin{Lm}\label{formula2}
With the notation above we have
$$
\Pi(s_{\th-\th_s})-\{\th-\th_s \}
\subseteq \{\a\in \RR_+\ |\ (\a^\vee,\th-\th_s)=1 \}\ .
$$
\end{Lm}
\begin{proof}
In the case when $p=3$, the matrix $\ring A$ equals $G_2$, and our statement
can be checked. Indeed, if $\a_1$, $\a_2$ is the standard basis, with
$\a_1$ the short root and $\a_2$ the long root, we have that
$\th=3\a_1+2\a_2$, $\th_s=2\a_1+\a_2$ and $\th-\th_s=\a_1+\a_2$.
Furthermore, $s_{\th-\th_s}=s_2s_1s_2$ is a reduced decomposition and
$$
\Pi(s_{\th-\th_s})=\{\a_1+\a_2,\ \a_2,\  3\a_1+2\a_2\}\ .
$$
Keeping in mind that
$$
(\a_1,\a_1)=2/3,\ \ \  (\a_1,\a_2)=-1\ \ \text{and} \ \  (\a_2,\a_2)=2,
$$
the conclusion follows.

When $p=2$, it is well known (see \cite[VI, \S 1, 3]{bourbaki})
that we have the following possible values for the scalar products
$(\a,\b^\vee)$ for any roots $\a,\b$ such that $\a\neq\pm \b$:
\begin{enumerate}
\item if $\b$ is long, $(\a,\b^\vee)\in\{0,\pm 1\}$;
\item if $\b$ and $\a$ are short, $(\a,\b^\vee)\in\{0,\pm 1\}$;
\item if $\b$ is short and $\a$ is long, $(\a,\b^\vee)\in\{0,\pm 2\}$.
\end{enumerate}
Moreover, $(\a,\a)\in\{1,2\}$.

Obviously, $\th-\th_s\in\Pi(s_{\th-\th_s})$. Let $\a\in\Pi(s_{\th-\th_s})$,
$\a\neq\th-\th_s$. Then,
$$
(s_{\th-\th_s}(\a),\th_s)=(\a,\th_s)  \ \ \text{and}\ \
(s_{\th-\th_s}(\a),\th)=(\a,2\th_s-\th)\ .
$$
Because $\a$ is a positive root, and $s_{\th-\th_s}(\a)$ is a negative root
the first scalar product must be zero. This implies that the second scalar
product equals $-(\a,\th)$. This one cannot be zero because it would
follow that $\a$ is fixed by $s_{\th-\th_s}$. Furthermore, $\a\neq \th$
because $(\th,\th_s)=1\neq 0$. Now, the above considerations on the
possible values of scalar products show that $(\a,\th)=1$. If we put all
these together, we get that
$$
(\a,(\th-\th_s)^\vee)=2\ .
$$
Keeping in mind that $\th-\th_s$ is a short root and it is different from
$\a$, the same considerations
on scalar products imply that $\a$ is a long root, and consequently
$(\a^\vee,\th-\th_s)=1$.
\end{proof}
We can choose a reduced decomposition for $s_{\th-\th_s}$
of the form
\begin{equation}\label{eq4}
s_{j_p}\cdots s_{j_1}s_{j_0}s_{j_1}\cdots s_{j_p},
\end{equation}
with $w=s_{j_1}\cdots s_{j_p}$ the minimal length element of
$\W$ for which $w(\th-\th_s)=\a_{j_0}$ is a simple (short) root.
Then, using formula (\ref{pi}) we see that
$$
\Pi(w)\subseteq \Pi(s_{\th-\th_s})-\{\th-\th_s \}.
$$
As before denote
$\Pi(w)=\{\a^{(k)}\ |\ 1\leq k\leq p\},$
with $\a^{(k)}=s_{j_p}\cdots s_{j_{k-1}}(\a_{j_k})$. By Lemma
\ref{formula2} we obtain that
$$
(\th-\th_s,(\a^{(k)})^\vee)=1,
$$
or equivalently
\begin{equation}\label{eq3}
(s_{j_{k+1}}\cdots s_{j_{p}}(\th-\th_s),\a_{j_k}^\vee)=1.
\end{equation}
\begin{Lm}\label{formula3}
With the notation above
$$
T_{w^{-1}}^{-1}X_{\th-\th_s}=X_{\a_{j_0}}T_w
$$
holds in the Artin group $\A_{\tilde W}$.
\end{Lm}
\begin{proof}
Using the formula (\ref{eq3}) and the relations in the double affine
Artin group we obtain that for all $p\geq k\geq 1$
$$
T_{\a_{j_k}}^{-1}X_{s_{j_{k+1}}\cdots s_{j_p}(\th-\th_s)}=
X_{s_{j_{k}}\cdots s_{j_p}(\th-\th_s)}T_{\a_{j_k}}.
$$
Now, our conclusion follows by applying these formulas.
\end{proof}
\begin{Prop}\label{formula4}
Let $\H_{\tilde W}$ be the double affine
group associated to an irreducible affine Cartan
matrix as before. Then, in $\H_{\tilde W}$ we have
\begin{equation}\label{eq9}
T^{-1}_{\th-\th_s}X_{\th-\th_s}-(T^{-1}_{\th-\th_s}X_{\th-\th_s})^{-1}=
t_s^{1/2}-t_s^{-1/2}.
\end{equation}
\end{Prop}
\begin{proof}
From (\ref{eq4}) we get that $T_{\th-\th_s}=T_{w^{-1}}T_{j_0}T_w$.
Therefore,
\begin{eqnarray*}
T^{-1}_{\th-\th_s}X_{\th-\th_s}&=&
T_{w}^{-1}T_{j_0}^{-1}T_{w^{-1}}^{-1}X_{\th-\th_s}\\
&=&T_{w}^{-1}T_{j_0}^{-1}X_{\a_{j_0}}T_w \\
&=&T_{w}^{-1}(X_{-\a_{j_0}}T_{j_0}+(t_s^{1/2}-t_s^{-1/2}))T_w\\
&=&X_{-\th+\th_s}T_{w^{-1}}T_{j_0}T_w+(t_s^{1/2}-t_s^{-1/2})\\
&=&(T^{-1}_{\th-\th_s}X_{\th-\th_s})^{-1}+(t_s^{1/2}-t_s^{-1/2}).
\end{eqnarray*}
We used Lemma \ref{formula3} and the fact that $\a_{j_0}$ is short.
\end{proof}

\medskip
\subsection{Proof of Theorem \ref{teorema2}}
Because we already have an isomorphism at the level of Artin groups all
we need to prove is that the relations (\ref{def1}) (\ref{def1.2}) 
(\ref{def1.3}) and (\ref{eq7})
are satisfied. In order to keep the computations in
$\H_{\tilde W}$ we will check the relations for
$\phi_{\tilde W^{\iota}}(T_j^\iota)$,
$\phi_{\tilde W^{\iota}}(T_j^\iota Y^\iota_{A_j^\iota})$ and 
$\phi_{\tilde W^{\iota}}(X^\iota_{-\a^\iota_j}T_j^\iota)$ 
(the corresponding relations for
$\phi_{\tilde W}(T_j)$ 
$\phi_{\tilde W}(T_jY_{A_j})$ and
$\phi_{\tilde W}(X_{-\a_j}T_j)$  
can be checked in the same way).
For $j\neq 0$ this is straightforward. For $j=0$
we have to show that
\begin{equation}\label{eq6}
\phi_{\tilde W^{\iota}}(X^\iota_{-\a_0^\iota}T_0^\iota)-
(\phi_{\tilde W^{\iota}}(X^\iota_{-\a_0^\iota}T_0^\iota))^{-1}=
(\phi_{\tilde W^{\iota}}(t_0))^{-1/2}-(\phi_{\tilde W^{\iota}}(t_0))^{1/2}.
\end{equation}
There are two possible situations:

{\sl Case 1.} If $p=1$, then $\a_j^\iota=\a_j$ for all $1\leq j\leq n$.
This implies that
\begin{eqnarray*}
\phi_{\tilde W^{\iota}}(X^\iota_{-\a_0^\iota}T_0^\iota)&=&
\phi_{\tilde W^{\iota}}
(qX^\iota_{\th^\iota}(T_{s_{\th^\iota}}^\iota)^{-1}Y^\iota_{-\th^\iota})\\
&=&q^{-1}Y_\th T_{s_\th}X_{-\th}\\
&=&T_0^{-1}X_{\a_0}.
\end{eqnarray*}
Therefore in this case the relation (\ref{eq6}) becomes precisely
the relation (\ref{eq7}).

{\sl Case 2.} If $p\neq 1$, then $\a_j^\iota=\frac{\a_j^\vee}{\sqrt{p}}$
for all $1\leq j\leq n$. As a consequence, with the notation in Section
\ref{combinatorics} we have that
$$
\th^\iota=\frac{\th_s^\vee}{\sqrt{p}}=\sqrt{p}\ \th_s\ \ \text{and}\ \
\th_s^\iota=\frac{\th}{\sqrt{p}}\ .
$$
Also, it is clear that $s_{\th^\iota}=s_{\th_s}$ and
$\phi_{\tilde W^{\iota}}(t_0)=t_s$. As before
\begin{eqnarray*}
\phi_{\tilde W^{\iota}}(X^\iota_{-\a_0^\iota}T_0^\iota)&=&
\phi_{\tilde W^{\iota}}
(qX^\iota_{\th^\iota}(T_{s_{\th^\iota}}^\iota)^{-1}Y^\iota_{-\th^\iota})\\
&=&q^{-1}Y_{\th_s^\vee} T_{s_{\th_s}}X_{-\th_s}.\\
\end{eqnarray*}
Using Lemma \ref{formula1} and the commuting relations in the double affine
Artin group we obtain
\begin{equation}\label{eq8}
\phi_{\tilde W^{\iota}}(X^\iota_{-\a_0^\iota}T_0^\iota)=
T_0^{-1}T_{s_{\th-\th_s}}^{-1}X_{\th-\th_s}T_0.
\end{equation}
Conjugating (\ref{eq6}) by $T_0$ we obtain the relation (\ref{eq9})
which was proved to be true in Proposition \ref{formula4}. The proof of
the Theorem \ref{teorema2} is completed.



\begin{thebibliography}{10}

\bibitem[B]{bourbaki}
{\sc N. Bourbaki}, {Groupes et alg\` ebres de Lie, Ch. IV, V, VI, }
{\sl Hermann}, Paris, 1968.


\bibitem[Br]{briesk}
{\sc E. Brieskorn}, Die Fundamentalgruppe des Raumes der regul\" aren
Orbits einer endlichen komplexen Spiegelungsgruppe,
{\sl Invent. Math.} {\bf 12} (1971), 57--61.



\bibitem[C1]{cherednik}
{\sc I. Cherednik}, Double affine Hecke algebras,
Knizhnik-Zamolodchikov equations, and Macdonald's operators,
{\sl Internat. Math. Res. Notices} 1992, no. 9, 171--180.

\bibitem[C2]{cherednik2}
{\sc I. Cherednik}, 
Double affine Hecke algebras and Macdonald's conjectures, 
{\sl Ann. of Math. (2)} {\bf 141} (1995), no. 1, 191--216. 


\bibitem[C3]{cherednik3}
{\sc I. Cherednik},
Intertwining operators of double affine Hecke algebras,
{\sl Selecta Math. (N.S.)} {\bf 3} (1997), no. 4, 459--495.


\bibitem[D]{deligne}
{\sc  P. Deligne},
{Les immeubles des groupes de tresses g\' en\' eralis\' es},
{\sl Invent. Math.} {\bf 17} (1972), 273--302.


\bibitem[H]{humph}
{\sc J. E. Humphreys}, {Reflection groups and Coxeter groups, }
{\sl Cambridge University Press}, Cambridge, 1990.


\bibitem[K]{kac}
{\sc V. G. Kac}, {Infinite dimensional Lie algebras (${\rm 3^{rd}}$ edition), }
{\sl Cambridge University Press},
Cambridge, 1990.

\bibitem[Kn]{knop}
{\sc F. Knop}, {Integrality of two variable Kostka functions, }
{\sl J. Reine Angew. Math.} {\bf 482} (1997),
177-189.

\bibitem[KS]{ks}
{\sc F. Knop and S. Sahi}, {A recursion and a combinatorial formula for Jack
polynomials,} {\sl Invent. Math.} {\bf 128} (1997) no.1, 9-22.

\bibitem[L1]{lek}
{\sc H. van der Lek},
Extended Artin groups,
{\sl Singularities, Part 2 (Arcata, Calif., 1981)}, 117--121,
Proc. Sympos. Pure Math., 40,
Amer. Math. Soc., Providence, RI, 1983.



\bibitem[L2]{lek-thesis}
{\sc H. van der Lek},
The homotopy type of complex hyperplane complements, Ph.D. Thesis,
Katholieke Universiteit Nijmegen, 1983.


\bibitem[Lo]{looijenga}
{\sc E. Looijenga},
Rational surfaces with an anti-canonical cycle,
{\sl Ann. of Math. (2)} {\bf 114} (1981), no. 2, 267--322.


\bibitem[Lu]{lusztig}
{\sc G. Lusztig}, Affine Hecke  algebras and their graded version,
{\sl J. Amer. Math. Soc.} {\bf 2} (1989), no. 3, 599--635.


\bibitem[M]{macdonald}
{\sc I. G. Macdonald}, Affine Hecke algebras and orthogonal polynomials, 
S\' eminaire Bourbaki, Vol. 1994/95, {\sl Ast\' erisque} 
{\bf 237} (1996), Exp. No. 797, 4, 189-207. 







\bibitem[S1]{s1}
{\sc S. Sahi}, {Interpolation, integrality and a generalization of Macdonald's 
polynomials,} {\sl Internat.
Math. Res. Notices} 1996, no. 10, 457-471.



\bibitem[S2]{sahi}
{\sc S. Sahi}, Nonsymmetric Koornwinder polynomials and duality,
{\sl Ann. of Math. (2)} {\bf 150} (1999), no. 1, 267--282.


\bibitem[ST]{saito2}
{\sc K. Saito and T. Takebayashi},
Extended affine root systems III. Elliptic Weyl groups,
{\sl Publ. Res. Inst. Math. Sci.} {\bf 33} (1997), no. 2, 301--329.





\end{thebibliography}
\end{document}